\documentclass[12pt,reqno]{amsart}
\usepackage[mathscr]{eucal}

\newcommand{\supp}{\mathop{\rm supp\,}}

\newcommand{\dist}{\mathop{\rm dist\,}}
\newcommand{\card}{\mathop{\rm card\,}}

\newcommand{\Int}{\mathop{\rm Int\,}}
\newcommand{\lip}{\mathop{\rm Lip\,}}
\newcommand{\N}{{\mathbb N}}

\newcommand{\C}{{C\{M_n\}}}
\newcommand{\D}{{\mathcal D}}
\newcommand{\E}{{\mathcal E}}

\newcommand{\al}{{\alpha}}

\newcommand{\e}{{\varepsilon}}
\renewcommand{\phi}{\varphi}

\newcommand{\rrr}[1]{{\rm (\ref{#1})}}

\newcommand{\xx}[2]{{x^{\langle #1\rangle}_{#2}}}
\newcommand{\xxx}[3]{{x^{\langle #1\rangle,#3}_{#2}}}

\newcommand{\spc}[1]{{{\rm Diff}^{#1}_0([0,1])}}

\newtheorem{thm}{Theorem}
\newtheorem{lem}{Lemma}
\newtheorem*{thma}{The growth gap theorem}

\theoremstyle{remark}
\newtheorem*{que1}{Question {\rm 1}}
\newtheorem*{que2}{Question {\rm 2}}
\newtheorem*{rem}{Remark}

\numberwithin{equation}{section}

\title[Distortion growth]
{Distortion growth for iterations of diffeomorphisms of the interval}
\author{Alexander Borichev}
\keywords{Iterations of diffeomorphisms, growth sequences.\newline 
$2000$ {\it Mathematical Subject Classification} 
37C05, 26A18, 58D05. 
\newline{\textit{\small E-mail}: \texttt{\small borichev@math.u-bordeaux.fr}}}
\begin{document}
\begin{abstract} 
We obtain several results on the distortion
asymptotics for the 
iterations of diffeomorphisms of the interval
extending the recent work of Polterovich and Sodin.
\end{abstract}
\maketitle

\section{\sc Main results}
\smallskip 

We consider the groups $\spc{N}$, $1\le N\le +\infty$, 
of all $C^N$-smooth diffeomorphisms
of the interval $[0,1]$ into itself fixing the endpoints $0,1$.
We associate to every $f\in \spc{1}$
its {\it growth sequence},
$$
\Gamma_n(f)=\max\bigl(\max_{x\in [0,1]}|(f^n)'(x)|,
\max_{x\in [0,1]}|(f^{-n})'(x)|\bigr), \qquad n\in\N,
$$
where $f^n$, $n\ge 1$, is the $n$-th iteration of $f$, 
and $f^{-n}$, $n\ge 1$, is the $n$-th iteration of the inverse 
diffeomorphism $f^{-1}$.
The asymptotics of the growth sequence does not change under conjugations:
for $g\in\spc{1}$,
$$
c(g)\Gamma_n(g^{-1}fg)\le \Gamma_n(f)\le C(g)\Gamma_n(g^{-1}fg), \qquad n\ge 1.
$$

This asymptotics is a basic dynamic invariant 
(see \cite{KH}). G.~D'Ambra and M.~Gromov proposed in \cite[7.10.C]{AG} 
to study the behavior of the growth sequences for 
various classes of diffeomorphisms
on smooth manifolds. Recently, L.~Polterovich and M.~Sodin \cite{PS} 
obtained several interesting results on
the growth sequences for diffeomorphisms in $\spc{2}$. 
In particular, they established 

\begin{thma}[{\cite[Theorem~1.7]{PS}}] If $f\in\spc{2}$, then 
either
$$
\lim_{n\to\infty}\frac{\log \Gamma_n(f)}n>0
$$
or
$$
\limsup_{n\to\infty}\frac {\Gamma_n(f)}{n^2}<\infty.
$$
\end{thma}

For other results on the growth sequences of diffeomorphisms
see \cite{P}, \cite{PS}, \cite{B}.

This paper is devoted to two problems related to the result 
of Poltero\-vich and Sodin. First, we would like to get more information
on the behavior of the growth sequences than that contained 
in the cited theorem,
possibly for smoother diffeomorphisms $f$, in terms of the 
local properties of $f$.
Another problem is to get analogs to the growth gap effect for diffeomorphisms 
of smoothness between $C^1$ and $C^2$. 

To formulate our results we need to introduce 
a decomposition of the set of fixed points
of a diffeomorphism. Let $1\le N\le \infty$, and let $f\in\spc{N}$. 
Denote by $E(f)$ the (closed) set of fixed points of $f$; 
$\{0,1\}\subset E(f)$. 
Consider the subsets of $E(f)$: 
\begin{align*}
E_1(f)&=\{x\in E(f):f'(x)\ne 1\},\\
E_k(f)&=\{x\in E(f):f'(x)=1,f''(x)=\ldots\\
&=f^{(k-1)}(x)=0,f^{(k)}(x)\ne 0\},
\qquad 1<k\le N.
\end{align*}
For $N=\infty$ we consider also
$$
E_\infty(f)=\{x\in E(f):f'(x)=1,f''(x)=\ldots=0\},
$$
and for $N<\infty$ we consider
$$
E^0_N(f)=\{x\in E(f):f'(x)=1,f''(x)=\ldots=f^{(N)}(x)=0\}.
$$
For $N=\infty$ we have
$$
E(f)=\bigsqcup_{1\le s\le\infty}E_s(f),
$$
and for $N<\infty$ we have
$$
E(f)=\bigsqcup_{1\le s\le N}E_s(f)\sqcup E^0_N(f),
$$
Set $V=\max_{x\in E(f)}|\log f'(x)|\ge 0$; $V>0$ if and only if
$E_1(f)\ne\emptyset$.

We use the notation $a(n)\sim b(n)$, $n\to\infty$, if 
$\lim_{n\to\infty}a(n)/b(n)=1$; $a(n)\asymp b(n)$, $n\to\infty$, if 
$0<c\le a(n)/b(n)\le C<\infty$.

Let $f\in\spc{1}$. It is known (see \cite[Appendix]{PS}) 
that if $E_1(f)\ne\emptyset$, then
$$
\log \Gamma_n(f)\sim n\,V,\qquad n\to\infty.
$$
Otherwise, if $E_1(f)=\emptyset$, then
$$
\log \Gamma_n(f)=o(n),\qquad n\to\infty.
$$

Furthermore (see \cite{PS}), for every $f\in\spc{1}$ different from
the identity map,
$$
\sum_{n\ge 1}\frac 1{\Gamma_n(f)}<\infty,
$$
and hence,
$$
\limsup_{n\to\infty}\frac{\log \Gamma_n(f)}{\log n}\ge 1.
$$
 
\begin{thm} Let $f\in\spc{\infty}$, and let $E_1(f)=\emptyset$. 
\begin{itemize}
\item[(A)] If $E_2(f)\ne\emptyset$, then
\begin{equation}
\Gamma_n(f)\asymp n^2,\qquad n\to\infty.\label{do1}
\end{equation}
\item[(B)] If $E_2(f)=\emptyset$, $E_\infty(f)\ne\emptyset$, then
\begin{equation}
\Gamma_n(f)=o(n^2),\qquad n\to\infty.\label{do2}
\end{equation}
\item[(C)] Finally, if $k\ge 3$, $E_s(f)=\emptyset$, $1<s<k$,
$E_\infty(f)=\emptyset$, $E_k\ne\emptyset$, then
$$
\Gamma_n(f)\asymp n^{k/(k-1)},\qquad n\to\infty.
$$
\end{itemize}
\label{t1}
\end{thm}

Of course, any $f\in\spc{\infty}$ with $E_1(f)=\emptyset$ satisfies 
one and only one condition among (A)--(C).

A version of this result for $f\in\spc{2}$ claims that 
if $E_1(f)=\emptyset$, $E_2(f)\ne\emptyset$, then \rrr{do1} holds,
and if $E_1(f)=E_2(f)=\emptyset$, then \rrr{do2} holds.
 
The principal part of Theorem~\ref{t1} is the part (B);
the parts (A) and (C) are rather standard; 
the lower estimates there follow, for example, from the description 
of the behavior of the iteration sequences $\{f^n(x)\}_{n\ge 1}$ 
near a fixed point of $f$ given in a book of Yu.~Lyubich
\cite[Section~2.6]{L}; the upper estimate in (A) follows from 
the growth gap theorem of \cite{PS}. 

A natural question is now whether additional smoothness conditions 
on $f\in\spc{\infty}$ may permit us to improve the asymptotics in the part 
(B) of Theorem~\ref{t1}. To answer this question in the negative, 
we fix a sequence $\{\e_k\}$ of positive numbers tending to $0$,
and a non-quasianalytic Carleman class $\C$,
$$
\C=\{f\in C^\infty([0,1]):
|f^{(n)}(x)|\le C(f)^nM_n,\,x\in[0,1]\},
$$
where $f^{(n)}$ is the $n$-th 
derivative of $f$, and $\log M_n$ is an increasing convex sequence,
$\lim_{n\to\infty}M_n=+\infty$. Recall that such a class is non-quasianalytic 
if for every two closed subintervals $I,J$ of $[0,1]$ with $I\subset \Int J$,
there exists $f\in\C$ with $0\le f\le 1$, $f{\bigm |}I\equiv 1$, and 
$\supp f\subset J$.
The Denjoy--Carleman theorem (see, for example, \cite{K,M}) claims that the
Carleman class $\C$ is non-quasianalytic if and only if 
$$
\sum_{n\ge 1}M_n^{-1/n}<\infty.
$$

\begin{thm} There exists $f\in\spc{\infty}\cap \C$ 
such that $E(f)=E_\infty(f)=\{0,1\}$, and
$$
\Gamma_n(f)\ge \e_n n^2,\qquad n\ge 1.
$$
\label{t2}
\end{thm}

Thus, we can conclude that some predictions from the Outlook of \cite{PS}
are true as demonstrated by Theorem~\ref{t1}; nevertheless, the
``optimistic scenario'' from the Outlook is disproved by Theorem~\ref{t2}.

\begin{que1} Suppose that $f\in\spc{\infty}$,
and that $E(f)=E_\infty(f)=\{0,1\}$.
What additional conditions should one impose on $f$ to
guarantee that
\begin{equation}
\limsup_{n\to\infty}\frac{\log \Gamma_n(f)}{\log n}=1?
\label{h1}
\end{equation}
\end{que1}

The previous theorem
shows that no additional smoothness conditions will work.
On the other hand, a bounded oscillation condition 
is sufficient for the property \rrr{h1} to hold. For simplicity,
we consider here only a model case.

\begin{thm} Let $f\in\spc{\infty}$ be different from
the identity map, and let $E(f)=E_\infty(f)=\{0,1\}$. 
Assume that $\phi(x)=f(x)-x>0$, $x\in(0,1)$. Suppose that for every $\e>0$
there exists $C_\e>0$ such that
\begin{equation}
\phi(y)\le C_\e(\phi(x))^{1-\e},\qquad 0<y<x<\frac12,
\label{h2}
\end{equation}
and the same inequality holds for $1/2<x<y<1$.
Then
$$
\limsup_{n\to\infty}\frac{\log \Gamma_n(f)}{\log n}=1.
$$
\label{t3n}
\end{thm}

Note that the condition \rrr{h2} guarantees, by a result of
F.~Sergeraert \cite{Se} that the germ of $f$ at $0$ 
imbeds in a flow of germs of $C^\infty$-smooth diffeomorphisms 
$f_\sigma$, $\sigma>0$, with $f_\sigma(x)-x$ flat at $0$
(and an analogous imbedding holds at the point $1$).
However, such an imbedding by itself does not provide, apparently, any new
information on the behavior of the growth sequence $\Gamma_n(f)$,
see the remark at the end of Section~\ref{s3}.

One more natural question here is whether there is the growth gap effect
for two (possibly non-commuting) diffeomorphisms 
$f,g$ in the group $\spc{\infty}$.
Denote by $[f,g]$ the subgroup of $\spc{\infty}$ generated by $f$ and $g$;
given $h\in[f,g]$ denote by $|h|_w$ the distance from the identity map
to $h$ in the word metric, that is the length of the shortest 
representation of $h$ by $f,f^{-1},g,g^{-1}$. 
Consider the growth sequence
$$
\Gamma_n(f,g)=\max_{h\in[f,g],\,|h|_w\le n}
\max_{x\in [0,1]}|h'(x)|,\qquad n\in\N.
$$

\begin{que2} Suppose that $\log \Gamma_n(f,g)=o(n)$, 
$n\to\infty$. Is it possible that the value
$$
\limsup_{n\to\infty}\frac{\log\log \Gamma_n(f,g)}{\log n}
$$
is positive? is equal to $1$?
\end{que2}

Next, we deal with diffeomorphisms of lower smoothness.
The class $\spc{1}$ is too large for any kind of the growth
gap effect to be present. Namely, given
a sequence $\{\e_k\}$ of positive numbers tending to $0$,
one can construct $f\in\spc{1}$ 
such that $E(f)=E^0_1(f)=\{0,1\}$, and
$$
\log \Gamma_n(f)\ge \e_n n,\qquad n\ge 1.
$$
 
The situation is different for the groupes 
$$
\spc{1,\al}=\spc{1}\cap C^{1,\al}, \qquad 0<\al\le 1,
$$
of diffeomorphisms with the derivative in the Lipschitz $\al$ class,
where
$$
C^{1,\al}=\bigl\{f\in C^1:|f'(x)-f'(y)|\le c(f)|x-y|^\al,\,
x,y\in[0,1]\bigr\}.
$$
(If $f\in\spc{1,\al}$, then automatically $f^{-1}\in C^{1,\al}$,
and hence $f^{-1}\in\spc{1,\al}$.)

For $\al=1$ we have just the growth gap of \cite{PS}:
if $f\in\spc{1,1}$, $E_1(f)=\emptyset$, then
$\Gamma_n(f)=O(n^2)$, $n\to\infty$.

In the case $0<\al<1$, we obtain a weaker form of the growth gap effect. 

\begin{thm} Let $0<\al<1$. {\rm (A)} 
If $f\in\spc{1,\alpha}$, and if $E_1(f)=\emptyset$, then
$$
\log\Gamma_n(f)=O(n^{1-\alpha}),\qquad n\to\infty.
$$

{\rm (B)} There exists $f_\alpha\in\spc{1,\alpha}$ such that
$E_1(f_\alpha)=\emptyset$, and
\begin{equation}
\lim_{n\to\infty}\frac{\log\log\Gamma_n(f_\alpha)}{\log n}=1-\alpha.\label{44}
\end{equation}
\label{t3}
\end{thm}

The plan of the paper is as follows.   
We start Section~\ref{s2} with an analysis of the behavior of the 
iteration sequences 
$\{f^n(x)\}_{n\ge 1}$ resembling the results contained 
in \cite[Section~2.6]{L}. 
In contrast to the asymptotical estimates of \cite{L} we obtain
global uniform estimates.
Then, using this analysis, 
we prove Theorems~\ref{t1} and \ref{t3n}.  
Theorem~\ref{t2} is proved in
Section~\ref{s3}. Finally, in Section~\ref{s4} we prove 
Theorem~\ref{t3}. Our proof
of the part (A) of Theorem~\ref{t3} imitates that of the original
growth gap theorem of Polterovich and Sodin \cite[Section~2]{PS}.
\bigskip

The author is thankful to Leonid Polterovich and Misha Sodin for numerous helpful discussions.

\section{\sc Proofs of Theorems~\ref{t1} and \ref{t3n}}
\label{s2}
\smallskip

{\bf 2.1.} First, we make several simple remarks.
Since $f$ is a diffeomorphism, we have
\begin{align*}
0<c(f)\le & f'(x)\le C(f),\\
-\infty<\log c(f)\le & \log f'(x)\le \log C(f),\qquad x\in[0,1].
\end{align*}
Let $x_1\in(0,1)$, $x_k=f(x_{k-1})$, $1<k\le n$.
We are going to estimate the value
$$
\Phi(n,x_1)=\Phi(n,x_1,f)=\sum_{k=1}^n\log f'(x_k),
$$
(and the supremum of $|\Phi|$ 
for $x_1\in(0,1)$) as a function of $n$.
This permits us to evaluate $a_n(f)$,
\begin{equation}
a_n(f)=\max_{[0,1]}\log[(f^n)'(x)]
\label{T}
\end{equation}
(and $a_n(f^{-1})$ is evaluated analogously).

Now we assume that $f$ is $C^2$-smooth, $E_1(f)=\emptyset$.
Replacing, if necessary, $f$ by $x\mapsto \Delta f(x/\Delta)$,
(and changing the domain $\D=[0,1]$ to $\D=[0,\Delta]$,
we can guarantee that
\begin{equation}
\max_\D|\phi''|\le 1.
\label{T1}
\end{equation}
Next we consider the set $\mathfrak A$ of
the closed subintervals $I$ of $\D$ such that
$E(f)\cap I=\partial I$. 
If $I\in\mathfrak A$, $\{x_k\}_{1\le k\le n}\cap I\ne\emptyset$, 
then $\{x_k\}_{1\le k\le n}\subset I$.
We have $f'=1$ at the endpoints of $I$. 

{\bf 2.2.} In the following three lemmas for simplicity of notations we assume
that the left end point of $I$ is $0$;
without loss of generality 
we assume that $f(x)\ge x$, $x\in I$. 
Set $\phi(x)=f(x)-x$, $x\in I=[0,b]$.
We have
$\phi(x)>0$ for $x\in (0,b)$,
$\phi(0)=\phi'(0)=\phi(b)=\phi'(b)=0$,
$\max_I\phi\le \Delta$, $\max_I|\phi''|\le 1$,
\begin{equation*}
-1<c(f)\le \min_I\phi'\le C(f).
\end{equation*}
Hence,
\begin{equation}
\bigl|\phi'(x)-\log(1+\phi'(x))\bigr|\le c(f)[\phi'(x)]^2,\qquad x\in I.
\label{do11}
\end{equation}
Furthermore, for $x,y\in I$,
\begin{multline}
\bigl|\log f(x)-\log f(y)\bigr|=
\Bigl|\log\frac{1+\phi'(x)}{1+\phi'(y)}\Bigr|\\
\le c(\phi)|\phi'(x)-\phi'(y)|\le c(\phi)|x-y|.\label{dododo}
\end{multline}
Since $\phi(x)=\phi'(x)=0$, $x\in\{0,b\}$, by \rrr{T1} we have
\begin{equation}
\phi(x)\le \min\{x^2/2,(b-x)^2/2\},\qquad x\in I.
\label{do6}
\end{equation}

\begin{lem} {\rm (A)} In this situation, we have
$$
\Bigl|\log\frac{\phi(x_n)}{\phi(x_1)}-\Phi(n-1,x_1)\Bigr|\le C(f)|I|,
$$
where $|I|$ is the length of $I$.
{\rm (B)} If $H>1$, and we have $\phi<1/100$, $H^{-1}-1<\phi'<H$ on $I$,
then we can choose $C(f)=C(H)$ in the previous inequality.
\label{pr1}
\end{lem}

Thus, to estimate $\Phi(n,x_1)$, we need only to know the asymptotics
of $\phi(x_n)$.

\begin{proof} (A) As a consequence of \rrr{do6} we have
$$
0\le x-(\phi(x))^{1/2}\le x+(\phi(x))^{1/2}\le b,\qquad x\in I.
$$

First we verify that for every $x\in I$, $\theta\in[0,1]$,
\begin{equation}
\bigl|\phi'(x+\theta(\phi(x))^{1/2})\bigr|\le 
3(\phi(x))^{1/2}.\label{1}
\end{equation}
Otherwise, (\ref{T1}) would imply that for some $x\in I$, 
the function $\phi'$ is of constant sign on the interval
$I_x=[x-(\phi(x))^{1/2},x+(\phi(x))^{1/2}]$, and
$$
|\phi'(y)|>(\phi(x))^{1/2},\qquad y\in I_x.
$$
If $\phi'(x)$ is positive, then $\phi(x-(\phi(x))^{1/2})<0$,
and if $\phi'(x)$ is negative, then $\phi(x+(\phi(x))^{1/2})<0$,
which is impossible. Thus, (\ref{1}) is proved.

Inequality (\ref{1}) implies that
\begin{equation}
|\phi'(x)|\le 3(\phi(x))^{1/2},\qquad x\in I,\label{2}
\end{equation}
and that
\begin{equation}
\frac{\phi(x+\theta(\phi(x))^{1/2})}{\phi(x)}\le 4, 
\qquad x\in I,\quad\theta\in[0,1].\label{3}
\end{equation}

Fix $x\in I$. Suppose that $\phi(x)<1/100$. 
For some $\theta,\theta_1\in[0,1]$,
\begin{gather*}
\int_{x}^{x+\phi(x)}\frac{\phi'(t)}{\phi(t)}dt=
\frac{\phi'(x+\theta\phi(x))}{\phi(x+\theta\phi(x))}\phi(x)=\\
\phi'(x+\theta\phi(x))\frac{\phi(x)}{\phi(x)+\theta\phi(x)
\phi'(x+\theta\theta_1\phi(x))}.
\end{gather*}
By (\ref{2}) and (\ref{3}) we obtain
$$
|\phi'(x+\theta\phi(x))|\le 
3(\phi(x+\theta\phi(x)))^{1/2}\le
6(\phi(x))^{1/2}\le \frac 35.
$$
Furthermore, 
\begin{multline*}
\Bigl|\frac{\phi(x)}{\phi(x)+\theta\phi(x)
\phi'(x+\theta\theta_1\phi(x))}-1\Bigr|=
\Bigl|\frac1{1+\theta\phi'(x+\theta\theta_1\phi(x))}-1\Bigr|\\
\le \frac 52\theta|\phi'(x+\theta\theta_1\phi(x))|
\le 15(\phi(x))^{1/2}.
\end{multline*}
Hence,
$$
\biggl|\int_{x}^{x+\phi(x)}\frac{\phi'(t)}{\phi(t)}dt-
\phi'(x+\theta\phi(x))\biggr|\le 90\phi(x),
$$
and for some $\theta_1\in[0,1]$,
\begin{multline}
\biggl|\int_{x}^{x+\phi(x)}\frac{\phi'(t)}{\phi(t)}dt-\phi'(x)\biggr|\\ \le 
90\phi(x)+|\phi''(x+\theta\theta_1\phi(x))|\cdot\theta\phi(x)\le
91\phi(x).\label{do4}
\end{multline}

Suppose now that $x\in\E=\{y\in\D:|\phi(y)|\ge 1/100\}$. The set $\E$
is a compact subset of $\D$ disjoint with $E(f)$. Furthermore,
for $y\in\E$, the interval $[y,y+\phi(y)]$ does not intersect $E(f)$. Hence, 
the function $y\mapsto \phi(y+\phi(y))$ does not vanish on $\E$, and as a consequence,
$$
|\phi(y+\phi(y))|\ge A(f)>0,\qquad y\in \E.
$$ 
Therefore,
\begin{gather}
\biggl|\int_{x}^{x+\phi(x)}\frac{\phi'(t)}{\phi(t)}dt-\phi'(x)\biggr|\le 
\Bigl|\log\frac{\phi(x+\phi(x))}{\phi(x)}\Bigr|+|\phi'(x)|\notag\\ \le 
\max\{\log\max_\D|100\phi|,\log\max_\D|\phi/A(f)|\}+\max_\D|\phi'|
\le c(f).\label{do5}
\end{gather}

Next, for $x=x_1,\ldots,x_{n-1}$ we sum up the inequalities
\rrr{do4} (if $x_k\not\in\E$) or \rrr{do5} (if $x_k\in\E$).
Since 
$$
\sum_{k\ge 1}\phi(x_k)\le |I|,
$$
we have 
$$
\card\bigl\{\{x_k\}_{1\le k\le n}\cap \E\bigr\}\le 100|I|.
$$
Using \rrr{do11} we obtain
\begin{gather*}
\Bigl|\log\frac{\phi(x_n)}{\phi(x_1)}-\Phi(n-1,x_1)\Bigr|\le
\sum_{k=1}^{n-1}\biggl|\int_{x_k}^{x_k+\phi(x_k)}\frac{\phi'(t)}{\phi(t)}dt-
\log(1+\phi'(x_k))\biggr|\\
\le c(f)\sum_{k=1}^{n-1}|\phi'(x_k)|^2+c(f)\sum_{k=1}^{n-1}\phi(x_k)+c(f)|I|.
\end{gather*}
Again using (\ref{2}), we conclude that
\begin{equation}
\Bigl|\log\frac{\phi(x_n)}{\phi(x_1)}-\Phi(n-1,x_1)\Bigr|\le C(f)|I|.
\label{do21}
\end{equation}
\smallskip

(B) In this case, $I\cap\E=\emptyset$, $c(f)=c(H)$ in \rrr{do11}, 
and the above argument shows that $C(f)=C(H)$ in \rrr{do21}.
\end{proof}

If we have more information on the size of $\phi''$, then we can extend
the estimate (\ref{3}) to bigger intervals.

\begin{lem} Let $0<\delta\le 1$. Fix $x\in I$.
Suppose that either {\rm (A)} $x\le b/2$ and 
\begin{equation}
\max_{[0,2x]}|\phi''(y)|\le \delta,
\label{do7}
\end{equation}
or {\rm (B)} $x>b/2$ and 
$$
\max_{[2x-b,b]}|\phi''(y)|\le \delta.
$$
Then $I_x=[x-\delta^{-1/2}(\phi(x))^{1/2},x+\delta^{-1/2}
(\phi(x))^{1/2}]\subset I$,
and
$$
\phi(y)\le 4\phi(x), \qquad y\in [x,x+\delta^{-1/2}(\phi(x))^{1/2}].
$$
\label{pr2}
\end{lem}

\begin{proof} The arguments in the cases (A) and (B) are 
analogous, and we restrict
ourselves to the case (A). By \rrr{do7},
$$
\phi(y)\le \delta y^2/2,\qquad y\in I,
$$
and hence, $I_x\subset [0,2x]\subset I$.
Let us verify that
\begin{equation}
\phi'(t)\le 3\delta^{1/2}
(\phi(x))^{1/2},\qquad t\in I_x.\label{1b}
\end{equation}
Otherwise, we would obtain that
$$
\phi'(y)>\delta^{1/2}(\phi(x))^{1/2},\qquad y\in I_x,
$$
and hence $\phi(x-\delta^{-1/2}(\phi(x))^{1/2})<0$, which is impossible.

Inequality \rrr{1b} implies our assertion: 
\begin{multline*}
\phi(y)-\phi(x)\le (y-x)\max_{x\le t\le y}\phi'(t)  \le
\delta^{-1/2}(\phi(x))^{1/2}\cdot 3\delta^{1/2}(\phi(x))^{1/2}\\=
3\phi(x),\qquad x\le y\le x+\delta^{-1/2}(\phi(x))^{1/2}.
\end{multline*}
\end{proof}

The following lemma shows that $x_n$ as a function of $n$
behaves asymptotically as the function inverse to the integral of
$1/\phi$. For related results see 
\cite[Appendix to Chapter~2, Theorem~3]{L}.

\begin{lem} Suppose that
\begin{equation}
\max_{[x_1,x_n]}|\phi'(x)|\le \frac12.
\label{s1}
\end{equation}
Then
$$
\frac 23\le
\frac1{n-1}\int_{x_1}^{x_n}
\frac{dt}{\phi(t)}\le 2.
$$
\label{pr3}
\end{lem}

\begin{proof} For every $x\in[x_1,x_{n-1}]$, $\theta\in[0,1]$,
there exists $\theta_1\in[0,1]$ such that
$$
\Bigl|\frac{\phi(x+\theta\phi(x))}{\phi(x)}-1\Bigr|=
\theta\bigl|\phi'(x+\theta\theta_1\phi(x))\bigr|\le \frac12.
$$
Since $x_{k+1}=x_k+\phi(x_k)$, we get
$$
\frac 23\le\int_{x_k}^{x_{k+1}}\frac{dt}{\phi(t)}\le 2.
$$
Summing up for $k=1,\ldots,n-1$ we obtain the assertion.
\end{proof}

{\bf 2.3.} We return to the analysis of the behavior of our sequence 
$\{x_k\}_{1\le k\le n}$
on the interval $I=[a,b]\in\mathfrak A$.
If $|f'-1|<1/2$ on $I$, then we put $J(I)=\emptyset$.
Otherwise, we choose the minimal  
closed subinterval $J(I)=[c,d]\subset I$
such that 
$$
|f'(x)-1|<\frac 12,\quad x\in I\setminus J(I)=I^l\cup I^r,
$$ 
where 
$$
I^l=[a,c),\qquad I^r=(d,b].
$$
Since $E_1(f)=\emptyset$, the set 
$$
J=\bigcup_{I\in\mathfrak A}J(I)
$$
is a compact subset of $\D$ having empty intersection with
$E(f)$.
Indeed, if $y\in\partial J$, then there are $y_k\in\partial J(I_k)$,
$I_k\in\mathfrak A$, such that $y_k\to y$ as $k\to\infty$.
Hence, $|f'(y)-1|\ge 1/2$, and $y\notin E(f)$. Therefore, for some 
$I\in\mathfrak A$,
we have $y\in \Int I$, and finally, $y\in J(I)\subset J$.
By continuity of $f$, 
$$
\inf_{x\in J}|f(x)-x|=\rho(f)>0.
$$
Since $f(x)-x$ is of constant sign on $I$, and 
$|x_{k+1}-x_k|\ge \rho(f)$ for $x_k\in J(I)$,
our sequence $\{x_k\}_{1\le k\le n}$ may contain at 
most $N(f)=1+\Delta/\rho(f)$ points of $J(I)$. 

Furthermore, if $J(I)=\emptyset$, then
condition (\ref{s1}) holds for our sequence $\{x_k\}_{1\le k\le n}\subset I$.
Otherwise, we have one of the following four possibilities: 
\smallskip

(I) either 
$\{x_k\}_{1\le k\le n}\subset I^l$
or $\{x_k\}_{1\le k\le n}\subset I^r$, and 
condition (\ref{s1}) holds for the sequence $\{x_k\}_{1\le k\le n}$;
\smallskip

(II) either $\{x_k\}_{1\le k\le n}\subset I^l\cup J(I)$
or $\{x_k\}_{1\le k\le n}\subset J(I)\cup I^r$.
Then dropping
at most $N(f)$ points of the sequence $\{x_k\}$ we return to 
the situation in (I) without changing the asymptotics of $\Phi$; 
\smallskip

(III) $\{x_k\}_{1\le k\le n}\cap I^l\ne\emptyset$, 
$\{x_k\}_{1\le k\le n}\cap I^r\ne\emptyset$,
$\{x_k\}_{1\le k\le n}\cap J(I)\ne\emptyset$.
Once again we assume that $f(x)\ge x$ on $I$.
Applying Lemma~\ref{pr1}, we get
$$
\Bigl|\log\frac{\phi(x_n)}{\phi(x_k)}-
\Phi(n-k,x_k)\Bigr|\le c(f), \qquad 1\le k\le n.
$$
If $j$ is the minimal index such that $x_j\in J(I)$, then
$$
\log\frac{\phi(x_n)}{\phi(x_j)}\le \log\frac{\max_\D\phi}{\rho(f)}\le c(f),
$$
and we are able to drop all the points $x_k$, $j\le k\le n$, and return to
(I) without worsening the asymptotics of $\Phi$; 
\smallskip

(IV) $\{x_k\}_{1\le k\le n}\cap I^l\ne\emptyset$, 
$\{x_k\}_{1\le k\le n}\cap I^r\ne\emptyset$,
$\{x_k\}_{1\le k\le n}\cap J(I)=\emptyset$.
Replacing $x_1$ by a suitable $x'_1\in [x_1,x_2]$,
and defining $x'_{k+1}=f(x'_k)$, $k\ge 1$,  
we can guarantee that $\{x'_k\}_{1\le k\le n}\cap 
J(I)\ne\emptyset$. Now, the sequence $\{x'_k\}_{1\le k\le n}$
satisfies the conditions of (II) or (III), and
$$
\bigl|\Phi(n,x_1)-\Phi(n,x'_1)\bigr|\le c(f)|I|.
$$
(Here we use \rrr{dododo}.)
\smallskip

Thus, from now on we may assume that the assertion of Lemma~\ref{pr3} 
holds for $\{x_k\}_{1\le k\le n}$.

\begin{proof}[{\bf 2.4.} Proof of Theorem~\rm\ref{t1}] The parts
(A) and (C) are rather standard; the lower estimates follow, for
example, from the description of the behavior of the iteration
sequences $\{f^n(x)\}_{n\ge 1}$ near a fixed point of $f$ given
in \cite[Section~2.6]{L}; to get the upper estimates
we can use an argument similar to that in the part (B).
\smallskip

(B) Fix $A\ge 1$. Let $\{x_k\}_{1\le k\le n}\subset
I=[a,b]\in\mathfrak A$, suppose that $\phi\ge 0$ on $I$, and let
$\{x_k\}_{1\le k\le n}$ and $I$ satisfy the conditions of
Lemma~\ref{pr3}. Using Lemma~\ref{pr2} and the fact that
$\phi''$ vanishes on $\partial I$ and is uniformly continuous on
$\D$, we obtain that for $\dist(x,\partial I)<\e(A,f)$,
\begin{equation}
\phi(t)\le 4\phi(x),\qquad x\le t\le x+A(\phi(x))^{1/2}.
\label{3d}
\end{equation}
If $\dist(x_1,\partial I)\ge\e(A,f)$,
then $\phi(x_1)\ge \beta(A,f)>0$, and then, by Lemma~\ref{pr1},
$$
\Phi(n-1,x_1)\le c(A,f).
$$
Otherwise, (\ref{3d}) holds for $x=x_1$.
Now, if $x_n\le x_1+A(\phi(x_1))^{1/2}$, then
$\phi(x_n)\le 4\phi(x_1)$, and again by Lemma~\ref{pr1},
$$
\Phi(n-1,x_1)\le c(f).
$$
Finally, if $x_n>x_1+A(\phi(x_1))^{1/2}$,
then by Lemma~\ref{pr3} and by (\ref{3d}),
$$
2n\ge \int_{x_1}^{x_n}\frac{dt}{\phi(t)}\ge
\int_{x_1}^{x_1+A(\phi(x_1))^{1/2}}\frac{dt}{\phi(t)}
\ge\frac A{4(\phi(x_1))^{1/2}}.
$$
Hence,
$$
\log\frac{\phi(x_n)}{\phi(x_1)}-2\log n\le
\log\Bigl(\frac{64\phi(x_n)}{A^2}\Bigr)\le
\log\Bigl(\frac{64}{A^2}\max_I\phi\Bigr).
$$
It remains to apply once again Lemma~\ref{pr1}
to conclude that 
$$
\Phi(n-1,x_1)\le 2\log n+c(f)-2\log A.
$$
Since $A$ is arbitrary, our proof is completed.
\end{proof}
 
\begin{proof}[{\bf 2.5.} Proof of Theorem~\rm\ref{t3n}] 
We argue as in the previous proof. Instead of Lemma~\ref{pr2}, we use 
the following resultat:

\begin{lem} If $\phi$ vanishes at $0$ with all its derivatives, 
$N\ge 1$, $x>0$, and
\begin{equation}
\phi(y)\le (\phi(x))^{1-1/(2N)}\le 1,\qquad 0<y<x,
\label{h4}
\end{equation}
then 
\begin{equation}
\phi(t)\le C(N,\phi)\phi(x),\qquad x\le t\le x+(\phi(x))^{1/N}.
\label{h3}
\end{equation}
\label{prn}
\end{lem}

Thus, condition \rrr{h4} permits us to extend the estimate
of Lemma~\ref{pr2} to much bigger intervals.

\begin{proof} The Gorny-Cartan inequalities (see, for example, 
\cite[6.4.IV]{M}) claim that for every $1\le k\le M$, and
$F\in C^M[0,1]$, there exists $C(M)$ such that
\begin{multline*}
|F^{(k)}(t)|\le C(M)\max_{[0,1]}|F(s)|^{1-k/M}\\
\times \max\Bigl[\max_{[0,1]}|F(s)|,\max_{[0,1]}|F^{(M)}(s)|
\Bigr]^{k/M},\qquad 0\le t\le 1,
\end{multline*}
where $F^{(k)}$ is the $k$-th 
derivative of $F$.

Applying these inequalities to $F(t)=\phi(xt)$, $M=4N^2$, $1\le k\le N$,
and using \rrr{h4} and the fact that for some $C(N,\phi)$,
\begin{align*}
|\phi^{(4N^2)}(t)|&\le C(N,\phi),\qquad 0\le t\le 1,\\
0\le \phi(t)&\le C(N,\phi)t^{4N^2},\qquad 0\le t\le 1,
\end{align*}
we conclude that
\begin{multline*}
\phi^{(k)}(x)=x^{-k}|F^{(k)}(1)|
\le C(N,\phi)x^{-k}(\phi(x))^{(1-1/(2N))(1-N/(4N^2))}\\
\le C_1(N,\phi)\cdot(\phi(x))^{1-1/N},\qquad 1\le k\le N,
\end{multline*}
for some $C_1(N,\phi)\ge 1$.

If $\phi(t)>C\phi(x)$ for some $x\le t\le x+(\phi(x))^{1/N}$, then
by induction we get a sequence of points $x\le t_{k+1}\le t_k\le t_0=t$, 
$1\le k< N$, with 
$$
\phi^{(k)}(t_k)>\bigl[C-kC_1(N,\phi)\bigr]
(\phi(x))^{1-k/N},\qquad 1\le k\le N.
$$
Fix $C(N,\phi)=NC_1(N,\phi)+\max_\D|\phi^{(N)}|$. Then we get
$$
\phi^{(N)}(t_N)>\max_\D|\phi^{(N)}|,
$$
which is impossible.
This contradiction proves \rrr{h3} with our choice of $C(N,\phi)$.
\end{proof}

Now we fix $N$ and $x_1\in\D$, 
and obtain as in the proof of Theorem~\ref{t1}~(B) that
either
$$
\Phi(n-1,x_1)\le c(N,f)
$$
or
$$
2n\ge \int_{x_1}^{x_n}\frac{dt}{\phi(t)}\ge
\int_{x_1}^{x_1+(\phi(x_1))^{1/N}}\frac{dt}{\phi(t)}
\ge\frac 1{C(N,\phi)(\phi(x_1))^{1-1/N}}.
$$
In the latter case,
$$
\Phi(n-1,x_1)\le \log\frac{\phi(x_n)}{\phi(x_1)}+c(\phi)\le
\frac{N}{N-1}\log n+c_1(N,\phi).
$$
Since $N$ is arbitrary, our proof is completed.
\end{proof}

\section{\sc Proof of Theorem~\ref{t2}}
\label{s3}
\smallskip
 
Without loss of generality we assume that $\{\e_k\}$
is a decreasing sequence.
We are going to construct 
a function $f$, $f(x)=x+\phi(x)$ with non-negative
$\phi\in\C$ vanishing at $0$
with all its derivatives, and points $\xx{n}{1}\in[0,1]$ such that
\begin{equation}
\Phi(n,\xx{n}{1})\ge \log (\e_n n^2),\qquad n\ge 1.
\label{4}
\end{equation}

The function $\phi$ will be of the form 
$t\mapsto \gamma_k(t-z_k)^2+\omega_k$ on disjoint intervals $J_k$
tending to $0$ with $\omega_k\ll \gamma_k\to 0$. In this way, we can make
$\max_{x\in[0,1]}\Phi(n,x)$ grow almost as fast as for $\phi(x)=x^2$,
and still keep $\phi$ smooth and flat at $0$.

We choose closed intervals $J_k$, and numbers $z_k,w_k$ such that
$$
[z_k,z_k+w_k]\subset\Int J_k\subset J_k\subset I_k=
\Bigl(\frac 1{2k},\frac 1{2k-1}\Bigr).
$$
Then we find non-negative $u_k\in\C$, $\supp u_k\subset I_k$, with 
$u_k{\bigm|}J_k\equiv 1$. Choose $1/100>\gamma_k\searrow 0$, $k\to\infty$, such that
for every sequence $\{\theta_k\}$, $\theta_k\in[0,1]$, the sums
$$
\sum_{k\ge 1}\gamma_k(\cdot-z_k)^2u_k
$$
and
$$
\sum_{k\ge 1}\theta_k\gamma_k u_k
$$
belong to $\C$ and vanish at $0$ with all their derivatives. 
Finally, we choose $n_k$ such that
\begin{gather}
\gamma_kw_k\ge \e_{n_k},\label{43}\\
n_k\ge \frac3{\gamma_k w_k}.\label{88}
\end{gather}

Now we put $\phi_k(x)=\gamma_k(x-z_k)^2$, $f_k(x)=x+\phi_k(x)$.
The functions $f_k$ satisfy the conditions of Lemma~\ref{pr1}~(B) with $H=2$.
We start with $\xxx{k}{1}{m}=z_k+m$, $0<m\le w_k/3$, and continue by 
$\xxx{k}{s+1}{m}=\xxx{k}{s}{m}+\phi_k({\xxx{k}{s}{m}})$, until 
$\xxx{k}{N(m)+1}{m}>z_k+2w_k/3$. Since $\gamma_k<1/100$, we have
$$
\xxx{k}{N(m)+1}{m}\le z_k+\frac{2w_k}3+\frac 1{100}\Bigl(\frac{2w_k}3\Bigr)^2<z_k+w_k.
$$
Using Lemma~\ref{pr3} we obtain that 
\begin{multline}
\frac 23N(m)\le\int_{m}^{\xxx{k}{N(m)+1}{m}-z_k}\frac{dt}{\gamma_k t^2}\le
\int_{m}^{w_k}\frac{dt}{\gamma_k t^2}\\
=\frac 1{\gamma_k}\Bigl[\frac 1m-\frac 1{w_k}\Bigr]
\le \frac23 \frac 1{\gamma_km}.
\label{brr}
\end{multline}
In particular, by \rrr{88},
$$
N(w_k/3)\le \frac 3{\gamma_k w_k}\le n_k.
$$
By continuity of $\phi_k$, for every $n\ge n_k$ there exists 
$m=m(n)\in(0,w_k/3]$ such that $n=N(m(n))$.
Furthermore, by Lemma~\ref{pr1}~(B),
\begin{gather*}
\Phi(n,\xxx{k}{1}{m(n)},f_k)\ge\log\frac{\gamma_k(\xxx{k}{n+1}{m(n)}-z_k)^2}
{\gamma_k(\xxx{k}{1}{m(n)}-z_k)^2}+C
\ge \log \frac{\gamma_k w_k^2}
{\gamma_k (m(n))^2}+C_1
\intertext{(by \rrr{43} and \rrr{brr})}
\ge\log [\e^2_{n_k}(N(m(n)))^2]+C_1=\log [\e^2_{n_k}n^2]+C_1,\quad n\ge n_k,
\end{gather*}
with $C$, $C_1$ independent of $k$ and $\{\e_n\}$.

Next we choose $0<\omega_k\le\gamma_k$ such that for $f^*_k$, 
$f^*_k(x)=\omega_k+\gamma_k(x-z_k)^2$, we still have
$$
\Phi(n,\xxx{k}{1}{m(n)},f^*_k)\ge \log (\e^2_{n_k}n^2)+C_1-1,
\qquad n_k\le n\le n_{k+1}. 
$$

It remains to define
$$
\phi_0(x)=\sum_{k\ge 1}\bigl[\gamma_k(x-z_k)^2+\omega_k\bigr]u_k(x),
$$
and add to $\phi_0$ a non-negative function in $\C$ vanishing at $0$
with all its derivatives, with support on 
$$
(0,1)\setminus\bigcup_{k\ge 1}[z_k,z_k+w_k],
$$
which is strictly positive on 
$$
(0,1)\setminus\bigcup_{k\ge 1}J_k,
$$
to get $\phi$. Now, (\ref{4}) is verified for all sufficiently big $n$.
Finally we change $\phi$ on a small interval inside $I_1$
to get (\ref{4}) for all $n\ge 1$.
 
\begin{rem} We can use the above construction to produce a flow $g^t$
of germs of $C^\infty$-smooth diffeomorphisms with $g^t(x)-x$ flat at $0$,
such that
\begin{equation}
\Gamma_n(g^1)\ge \e_n n^2,\qquad n\ge 1.
\label{h5}
\end{equation}

To do this, consider the equation
$$
\left\{
\begin{gathered}
\frac{\partial F}{\partial t}(t,x)=\phi(F(t,x)),\qquad x,t\ge 0,\\
F(0,x)=x,\qquad x\ge 0.
\end{gathered}
\right.
$$
Then $g^t=F(t,\cdot)$, $t\ge 0$, are the germs of 
$C^\infty$-smooth diffeomorphisms, and  the germs $g^t(x)-x$ are flat at $0$.
Put $g=g^1$. An easy argument shows that for $x>0$, $n\ge 1$, such 
that $g^n(x)$ is sufficiently small, we have
$$
\int_x^{g^n(x)}\frac{dt}{\phi(t)}=
\int_{F(0,x}^{F(n,x)}\frac{dt}{\phi(t)}=
\int_{0}^{n}\frac{\frac{\partial}{\partial s}F(s,x)}{\phi(F(s,x))}ds=n.
$$
and hence
$$
(g^n)'(x)=\frac{\phi(g^n(x))}{\phi(x)}.
$$
Starting from these equalities, and using the same argument 
as in the previous proof, we conclude that \rrr{h5} holds.
\end{rem}

\section{\sc Proof of Theorem~\ref{t3}}
\label{s4}
\smallskip

(A) We use the scheme proposed in \cite{PS}. First we prove two lemmas:
a Denjoy-type statement and a convex analysis result. 

\begin{lem} Let $f\in\spc{1,\alpha}$. If $J\subset [0,1]$ is a closed interval
such that $f(J)\cap J=\emptyset$, then for every $n\in\N$ and every $x,y\in J$,
\begin{equation}
\Bigl|\log\frac{(f^n)'(x)}{(f^n)'(y)}\Bigr|\le 
c(f)n^{1-\alpha}.
\label{42}
\end{equation}
\label{l4}
\end{lem}

\begin{proof} Since $0<\min_{[0,1]}f'\le \max_{[0,1]}f'<\infty$, 
$f\in C^{1,\alpha}$,
we have
\begin{equation}
\Bigl|\log\frac{(f^n)'(x)}{(f^n)'(y)}\Bigr|\le 
c(f)\sum_{k=1}^n |f'(x_k)-f'(y_k)|\le
c_1(f)\sum_{k=1}^n |x_k-y_k|^\alpha,
\notag
\end{equation}
where $x_1=x$, $y_1=y$, $x_k=f(x_{k-1})$, $y_k=f(y_{k-1})$, $k>1$. 
Next, the intervals $[x_k,y_k]$ are disjoint, and hence
$\sum_{k=1}^n |x_k-y_k|\le 1$.
By the H\"older inequality we obtain
$$
\sum_{k=1}^n |x_k-y_k|^\alpha \le \Bigl(\sum_{k=1}^n |x_k-y_k|\Bigr)^\alpha 
\Bigl(\sum_{k=1}^n 1\Bigr)^{1-\alpha}\le n^{1-\alpha}, 
$$
and (\ref{42}) follows.
\end{proof}
 
\begin{lem}[compare to Lemma~2.3 of \cite{PS}] 
Let $\{a_n\}_{n\ge 0}$ be a sequen\-ce of non-negative numbers.
Suppose that the following 
almost convexity inequality
\begin{equation}
2a_n-a_{n-1}-a_{n+1}\le K\exp\bigl[-a_n+K_1n^{1-\alpha}\bigr],\qquad n\ge 1,
\label{41}
\end{equation}
holds for some positive $K,K_1$, $a_0=0$, and 
\begin{equation}
\liminf_{n\to\infty}\frac{a_n}{n}=0.
\label{38}
\end{equation}
Then
\begin{equation}
a_n\le An^{1-\alpha},\qquad n\ge 1,
\label{do15}
\end{equation}
where $A=A(K,K_1)$.
\label{l5}
\end{lem}

\begin{proof} First we fix $A$ so large that for every $n\ge 1$,
\begin{multline}
A\bigl[(n+1)^{1-\alpha}-n^{1-\alpha}\bigr]\ge 
AC(\alpha)n^{-\alpha}\\
\ge
2K\sum_{k>n}\exp\Bigl[-\Bigl(\frac A2-K_1\Bigr)k^{1-\alpha}\Bigr].
\label{do12}
\end{multline}
If \rrr{do15} does not hold, then
we could find the smallest integer $n$ such that 
\begin{equation}
a_n\le An^{1-\alpha},\qquad a_{n+1}>A(n+1)^{1-\alpha},
\label{36}
\end{equation}
and hence,
\begin{equation}
a_{n+1}-a_n>A\bigl[(n+1)^{1-\alpha}-n^{1-\alpha}\bigr].
\label{do18}
\end{equation}

Then either
\begin{equation}
a_k\ge \frac A2k^{1-\alpha},\qquad k>n,
\label{do16}
\end{equation}
or we could find the smallest integer $m>n$ such that
\begin{equation}
a_m\ge \frac A2m^{1-\alpha},\qquad a_{m+1}<\frac A2(m+1)^{1-\alpha}.
\label{37}
\end{equation}
In the latter case, by (\ref{41}),
\begin{multline*}
(a_{n+s+1}-a_{n+s})-(a_{n+s+2}-a_{n+s+1})\le \\
K\exp\Bigl[-\Bigl(\frac A2-K_1\Bigr)(n+s+1)^{1-\alpha}\Bigr],
\qquad 0\le s<m-n.
\end{multline*}
Now, by \rrr{do12} and by \rrr{do18},
\begin{multline}
a_{n+s+2}-a_{n+s+1}\ge\\A\bigl[(n+1)^{1-\alpha}-n^{1-\alpha}\bigr]-
K\sum_{k>n}\exp\Bigl[-\Bigl(\frac A2-K_1\Bigr)k^{1-\alpha}\Bigr]
\ge\\ \frac A2\bigl[(n+1)^{1-\alpha}-n^{1-\alpha}\bigr],\qquad 0\le s<m-n,
\label{40}
\end{multline}
and
\begin{multline*}
a_{m+1}=a_{n+1}+\sum_{k=n+1}^{m}[a_{k+1}-a_{k}]\\ \ge
A(n+1)^{1-\alpha}+
(m-n)\frac A2\bigl[(n+1)^{1-\alpha}-n^{1-\alpha}\bigr].
\end{multline*}
Since the function $x\mapsto x^{1-\alpha}$ is concave, we obtain
$$
a_{m+1}\ge \frac A2(n+1)^{1-\alpha}+
\frac A2\sum_{k=n+1}^{m}[(k+1)^{1-\alpha}-k^{1-\alpha}]=
\frac A2(m+1)^{1-\alpha},
$$
and we get a contradiction to (\ref{37}). Thus, \rrr{do16} is established.
Arguing as above we derive from (\ref{40}) that
$$
\frac {a_k}{k}\ge \frac A2 \bigl[(n+1)^{1-\alpha}-n^{1-\alpha}\bigr]+o(1),
\qquad k\to\infty,
$$
that contradicts to our condition (\ref{38}). Thus, (\ref{36})
does not hold, and \rrr{do15} is proved.
\end{proof}

Now, the assertion (A) of Theorem~\ref{t3} follows just as in 
\cite{PS}. To make our
proof self-contained, we repeat here the argument from \cite{PS}. 

Consider the sequence $a_n(f)$ defined by \rrr{T}. Fix $n\ge 1$, and 
choose $x_1\in[0,1]$
such that for $x_{k+1}=f(x_k)$, $1\le k\le n$, we have
$$
a_n(f)=\log[(f^n)'(x_2)]=\sum_{k=2}^{n+1}\log f'(x_k).
$$
Then
\begin{align*}
a_{n+1}(f)&\ge\sum_{k=1}^{n+1}\log f'(x_k),\\
a_{n-1}(f)&\ge\sum_{k=3}^{n+1}\log f'(x_k),
\end{align*}
and as a result,
\begin{gather*}
2a_n(f)-a_{n+1}(f)-a_{n-1}(f)\le \log f'(x_2)-\log f'(x_1)\\
\le c(f)|x_2-x_1|^\alpha\le c(f)\Bigl|\frac{x_2-x_1}{x_{n+2}-x_{n+1}}
\Bigr|^\alpha=
\frac{c(f)}{|(f^n)'(y)|^\alpha}
\end{gather*}
for some $y$ between $x_1$ and $x_2$. By Lemma~\ref{l4},
$$
|(f^n)'(y)|\ge |(f^n)'(x_2)|\exp[-c(f)n^{1-\al}]=\exp [a_n(f)-c(f)n^{1-\al}],
$$
and we conclude that the sequence $a_n(f)$ satisfies the condition \rrr{41}.
Since $E_1(f)=\emptyset$, this sequence satisfies also \rrr{38}, and 
we can apply Lemma~\ref{l5} to complete the proof of the part (A) of 
the theorem. 
\smallskip

(B) First we fix $\beta>0$, and consider the function
$$
\bar\phi(x)=\bar\phi(\beta,x)=(x^{-1/\beta}-1)^{-\beta}-x-
x^{(\alpha+1)(\beta+1)/\beta}\sin\frac{2\pi}{x^{1/\beta}},\quad 0<x<1.
$$
Since $\alpha>0$, we have
$$
\bar\phi(x)\sim \beta x^{(\beta+1)/\beta},\qquad x\to 0.
$$

Fix a positive integer $N$. If $x_1=N^{-\beta}$, then
$$
x_2=x_1+\bar\phi(x_1)=(N-1)^{-\beta},
$$
and by induction we obtain
$$
x_k=x_{k-1}+\bar\phi(x_{k-1})=(N+1-k)^{-\beta},\qquad 1\le k\le N.
$$

Furthermore,
\begin{multline*}
\bar\phi'(x)=(1-x^{1/\beta})^{-\beta-1}-1
-(\alpha+1)\frac{\beta+1}\beta
x^{[(\alpha+1)(\beta+1)/\beta]-1}\sin\frac{2\pi}{x^{1/\beta}}\\
+\frac{2\pi}\beta x^{\alpha(\beta+1)/\beta}\cos\frac{2\pi}{x^{1/\beta}},
\qquad x>0,
\end{multline*}
and
$$
\bar\phi'(k^{-\beta})=\Bigl(1-\frac 1k\Bigr)^{-\beta-1}-1
+\frac{2\pi}\beta k^{-\alpha(\beta+1)},\qquad k\in \N.
$$

If $\beta+1<1/\alpha$, then
\begin{align}
\bar\phi'(k^{-\beta})\sim &\frac{2\pi}\beta k^{-\alpha(\beta+1)},
\qquad k\in \N,\quad k\to\infty,\label{741}\\
\bar\phi'\Bigl(\bigl(k+\frac 12\bigr)^{-\beta}\Bigr)\sim &
-\frac{2\pi}\beta k^{-\alpha(\beta+1)},\qquad k\in \N,
\quad k\to\infty.\notag
\end{align}
In particular, we obtain that $\bar\phi'$ vanishes on a sequence of
points $y_k$,
$$
\bigl(k+\frac 12\bigr)^{-\beta}<y_k<k^{-\beta},\qquad k>1.
$$

The formula \rrr{741} implies that
\begin{equation}
\sum_{k=1}^N \log(1+\bar\phi'(x_k))\sim \frac{2\pi}{\beta (1-\alpha(\beta+1))}
N^{1-\alpha(\beta+1)},
\quad N\to\infty.
\label{39}
\end{equation}

Next we verify that $\bar\phi'$ belongs to the Lipschitz $\alpha$ class.

\begin{lem} If $p>0$, $0<\alpha<1$, $b\ge (p+1)\alpha$, then the 
functions $g,h$,
$$
g(x)=x^b\sin x^{-p},\qquad h(x)=x^b\cos x^{-p},
$$
belong to $\lip_\alpha[0,1]$.
\label{l6}
\end{lem}

\begin{proof} Let $0\le y< x \le 1$.
If $0<x-y\le x^{p+1}$, then
$$
|g(x)-g(y)|\le (x-y)\cdot\max_{y\le t\le x}|g'(t)|\le c(x-y)x^{b-p-1}
\le c_1(x-y)^\alpha,
$$
where $c,c_1$ depend only on $b,p,\alpha$.
Otherwise, if $y<x-x^{p+1}$, then
$$
|g(x)-g(y)|\le 2\max_{y\le t\le x}|g(t)|\le 
2x^b\le 2(x-y)^\alpha.
$$
The same argument works for the function $h$.
\end{proof}

Since the function $x\mapsto (1-x^{1/\beta})^{-\beta-1}$ is
$C^1$-smooth on $[0,1/2]$, we conclude that
$$
\|\bar\phi'\|_{\lip_\alpha[0,1/2]}\le K(\beta).
$$

Finally, we use the functions $\bar\phi(\beta,\cdot)$, $\beta>0$,
to construct $f$ with
$$
\lim_{N\to\infty}N^{\alpha(\beta+1)-1}\max_{x\in[0,1]}\Phi(N,x,f)>0
$$
for every $\beta>0$.

Fix $\psi\in C^\infty([0,1])$, $0\le\psi\le 1$, with $\supp\psi\subset[0,1)$,
$\supp(1-\psi)\subset(0,1]$. Choose a sequence $\dfrac 1\al-1>\beta_k\to 0$,
and a sequence of disjoint intervals $I_k=[a_k,b_k]\subset [0,1]$, $k\ge 1$.
For every $k\ge 1$ we define $\phi(\cdot)=\bar\phi(\beta_k,\cdot)$,
$\Delta=|I_k|(K(\beta_k))^{-1/\al}$.
Given $0<\delta<\Delta/2$ such that $\phi'(\delta)=0$, define
$$
\phi_\delta(x)=
\begin{cases}
\phi(x),\qquad 0<x\le\delta,\\
\phi(\delta),\qquad \delta<x\le \Delta/2,\\
\phi(\delta)\cdot\psi((2x-\Delta)/\Delta),\qquad \Delta/2\le x\le \Delta.
\end{cases}
$$
Then 
$$
\|\phi'_\delta\|_{\lip_\alpha[0,\Delta]}\le K(\beta_k)+c|\phi(\delta)|
\Delta^{-1-\al},
$$
and for sufficiently small $\delta$ with $\phi'(\delta)=0$,
$$
\|\phi'_\delta\|_{\lip_\alpha[0,\Delta]}\le 2K(\beta_k).
$$
Fix such a value $\delta=\delta(k)$, and define
$$
\phi^*_k(x)=\frac{|I_k|}{\Delta}\phi_{\delta}\Bigl
(\frac{\Delta x}{|I_k|}\Bigr),\qquad 0\le x\le |I_k|.
$$
Then
$$
\|{\phi^*_k}'\|_{\lip_\alpha[0,|I_k|]}\le 2\frac{\Delta^\al}{|I_k|^\al}
K(\beta_k)\le 2.
$$
Put
$$
f(x)=\begin{cases}
x+\phi^*_k(x-a_k),\qquad x\in I_k,\quad k\ge 1,\\
x\qquad \text{elsewhere}.
\end{cases}
$$
Then
$$
\|f'\|_{\lip_\alpha[0,1]}\le 4,
$$
and for every $k\ge 1$, by \rrr{39} we obtain
$$
\Phi(N,a_k+N^{-\beta_k},f)\sim \frac{2\pi}{\beta_k(1-\alpha(\beta_k+1))}
N^{1-\alpha(\beta_k+1)},\qquad N\to\infty.
$$
Hence, for every $k\ge 1$,  
$$
\liminf_{N\to\infty}\frac{\log\log\Gamma_N(f)}{\log N}\ge 1-\alpha(\beta_k+1).
$$
Applying the result of the part (A) we obtain that $f$ 
satisfies (\ref{44}), and the proof of Theorem~\ref{t3} is completed.
 
\begin{rem} Note that $\phi$ and $\phi'$
constructed in the part (B) of the previous proof
may vanish at $0$ more rapidly than any preassigned power 
if we take sufficiently small $\beta>0$.
However, $\phi''$ is always unbounded near the point $0$.
\end{rem}

\bigskip

\noindent \textsf{\small Alexander Borichev, Department of Mathematics,} 
\newline
\noindent\textsf{\small University
of Bordeaux I, 351, cours de la Lib\'eration, 33405 Talence, France}

\end{document}